\documentclass[11pt]{article}
\usepackage{latexsym,amssymb}
\usepackage{amsmath}
\usepackage{color}

\newtheorem{Theorem}{\bf Theorem}[section]
\newtheorem{Lemma}{\bf Lemma}[section]
\newtheorem{Proposition}{\bf Proposition}[section]
\newtheorem{Corollary}{\bf Corollary}[section]
\newtheorem{Remark}{\bf Remark}[section]
\newtheorem{Example}{\bf Example}[section]
\newtheorem{Definition}{\bf Definition}[section]

\newenvironment{theorem}{\begin{Theorem}$\!\!\!$}{\end{Theorem}}
\newenvironment{lemma}{\begin{Lemma}$\!\!\!$}{\end{Lemma}}

\newenvironment{corollary}{\begin{Corollary}$\!\!\!$}{\end{Corollary}}
\newenvironment{remark}{\begin{Remark}$\!\!\!$}{\end{Remark}}

\newenvironment{definition}{\begin{Definition}$\!\!\!$}{\end{Definition}}

\numberwithin{equation}{section}

\begin{document}
\title{Critical exponent for the global existence of solutions\\
to a nonlinear degenerate/singular parabolic equation}
\author{
Yohei Fujishima\\
Department of Mathematical and Systems Engineering,\\ 
Shizuoka University\\
3-5-1 Johoku, Hamamatsu 432-8561, Japan
\\
\\
Tatsuki Kawakami\\
Department of Applied Mathematics and Informatics,\\
Ryukoku University\\
Seta Otsu 520-2194, Japan
\\
\\
Yannick Sire\\
Department of Mathematics,\\
Johns Hopkins University\\
3400 N. Charles Street
Baltimore, MD 21218,USA
}
\date{}
\maketitle

\begin{abstract}
We investigate a non-homogeneous nonlinear heat equation which involves degenerate 
or singular coefficients belonging to the $A_2$ class of functions. 
We prove the existence of a Fujita exponent 
and describe the dichotomy existence/non-existence of global in time solutions. 
The $A_2$ coefficient admits either a singularity at the origin or a line of singularities. 
In this latter case, the problem is related to the fractional laplacian, through the Caffarelli-Silvestre extension 
and is a first attempt to develop a parabolic theory in this setting.
\end{abstract}

\tableofcontents
%
\section{Introduction}
We consider the problem
\begin{equation}
\label{eq:1.1}
\left\{
\begin{array}{ll}
\displaystyle{\partial_t u-\mbox{div}(w(x)\nabla u)=u^p}, & x\in{\mathbb R}^N,\,\,\,t>0,\vspace{5pt}\\
\displaystyle{u(x,0)=\varphi(x)}\ge0,\qquad & x\in{\mathbb R}^N,
\end{array}
\right.
\end{equation}
where the coefficient $w$ is either $w(x)=|x_1|^a$
with $a\in(-1,1)$, or $w(x)=|x|^b$ with $b \in (-N,N)$. Here one has $N\ge1$,
$\partial_t:=\partial/\partial t$ and $p>1$. 

The aim of the present work is to develop a global-in-time existence theory of mild solutions 
for the problem \eqref{eq:1.1}. 
The coefficient $w(x)$ depending on the powers $a$ or $b$ degenerates or blows up. 
We prove that there is a critical exponent for the global existence of 
positive solutions of problem \eqref{eq:1.1}, the so-called Fujita exponent. 
\vspace{5pt}

We give first the definition of a solution to \eqref{eq:1.1}.
\begin{definition}
\label{Definition:1.1}
Let $\varphi$ be a nonnegative measurable function in ${\mathbb R}^N$.
Let $T\in(0,\infty]$ and $u$ be a nonnegative measurable function in ${\mathbb R}^N\times(0,T)$ such that
$u\in L^\infty(0,T: L^\infty({\mathbb R}^N))$.
Then we call $u$ a solution of \eqref{eq:1.1} in ${\mathbb R}^N\times(0,T)$
if $u$ satisfies
\begin{equation}
\label{eq:1.2}
u(x,t)=\int_{{\mathbb R}^N}\Gamma(x,y,t)\varphi(y)\,dy
+\int_0^t\int_{{\mathbb R}^N}\Gamma(x,y,t-s)u(y,s)^p\,dy\,ds<\infty
\end{equation}
for almost all $x\in{\mathbb R}^N$ and $t\in(0,T)$.
In particular, we call $u$ a global-in-time solution of \eqref{eq:1.1}
if $u$ is a solution of \eqref{eq:1.1} in ${\mathbb R}^N\times(0,\infty)$.
Here $\Gamma=\Gamma(x,y,t)$ is the fundamental solution of
$$
\partial_t v-{\rm{div}}(w(x)\nabla v)=0,\qquad x\in {\mathbb R}^N,\quad t>0,
$$
with pole at $(y,0)$.
\end{definition} 
\vspace{8pt}
The previous definition is the well-known class of {\sl mild solutions} and is natural to tackle parabolic problems. 
A main point of the previous definition is that 
it involves the fundamental solution of the operator under consideration. 
It is important to notice that in our context, due to the non-homogeneity of the operator, 
the fundamental solution is not translation-invariant. 
Furthermore, there is no explicit expression of it, though bounds are known. This makes the theory harder. 

We discuss now the features of the weight $w(x)$. In both cases under consideration, 
the weights belong to the class $A_2$ of Muckenhoupt functions \cite{muck}. 
This class of functions is very important in harmonic analysis for the boundedness of Maximal Functions. 
From the PDE point of view, elliptic equations and potential theory involving these weights 
have been studied in \cite{FKS1,FKS2,FKS3}. 
See also \cite{CS} for the parabolic counterpart. 

In the present work, we do not consider general weights 
since it is very complicated in this case to give precise results as our aim is. 
We will consider two types of weights.  
The first one is $|x_1|^a$ which is $A_2$ if and only if $a \in (-1,1)$. 
This exhibits singularities along the line $x_1=0$. 
The other weight under consideration is $|x|^b$ which is $A_2$ for $a \in (-N,N)$ 
and exhibits a singularity at the origin $x=0$. 
This former function is particularly interesting since this is related to equations involving fractional laplacians. 
Indeed, Caffarelli-Silvestre \cite{cafSil} proved that 
the fractional laplacian $(-\Delta)^s$ is the Dirichlet-to-Neumann operator of a suitable extension 
in the half-space involving the operator $y^{1-2s}$ which is degenerate along the line $y=0$. 
Notice that since $s \in (0,1)$ then $1-2s \in (-1,1)$. 
The Fujita problem for the fractional laplacian,
that is,
\begin{equation}
\label{eq:1.3}
\partial_tu+(-\Delta)^su=u^p,\qquad x\in {\mathbb R}^N,\quad t>0,
\qquad
u(x,0)=\varphi(x)\ge0,
\qquad x\in {\mathbb R}^N
\end{equation}
was studied in several papers.
Among others, Sugitani~\cite{S} showed the Fujita exponent for this problem (see also \cite{IKK01}).
In the problem~\eqref{eq:1.3},
the operator is non-local, but we can use an explicit form of the fundamental solution.
On the other hand, in the problem~\eqref{eq:1.1},
we can't use it even the operator is local,
and several fundamental topics have been left open up to now.
One of the goals of the present paper is an attempt to fully understand the parabolic theory 
for the fractional laplacian by considering these degenerate weights. 
\vspace{8pt}

We now describe our results. We first introduce some notations.
For any $x\in{\mathbb R}^N$ and $r>0$,
we put $B_r(x):=\{y\in{\mathbb R}^N : |x-y|<r\}$.
For any $1\le r\le\infty$, we denote by $\|\cdot\|_r$ the usual norm
of $L^r:=L^r({\mathbb R}^N)$.
For any measurable function $f$ in ${\mathbb R}^N$,
$$
\mu(\lambda):=\left|\{x\,:\,|f(x)|>\lambda\}\right|,
\qquad \lambda\ge 0,
$$
is the distribution function of $f$, and we define
the non-increasing rearrangement of $f$ by
$$
f^*(s):=\inf\{\lambda>0\,:\,\mu(\lambda)\le s\}.
$$
The spherical rearrangement of $f$ is defined by
$$
f^\sharp(x):=f^*(c_N|x|^N),
$$
where $c_N$ is the volume of the unit ball in ${\mathbb R}^N$.
For any $1\le r\le\infty$,
we define the $L^{r,\infty}$~space by
$$
L^{r,\infty}:=\{f\,:\, \mbox{$f$ is measurable in ${\mathbb R}^N$},\,\,\, \|f\|_{r,\infty}<\infty\},
$$
where
$$
\|f\|_{r,\infty}:=\sup_{s>0}\,s^{1/r}f^{**}(s),
\qquad
f^{**}(s):=\frac{1}{s}\int_0^s f^*(r)dr.
$$
Then $L^{r,\infty}$ is a Banach space and the following holds (see e.g. \cite{G}):
\begin{itemize}
  \item Let $1<r<\infty$. Then $f\in L^{r,\infty}$ if and only if
  $$
  0\le f^\sharp(x)\le C_1|x|^{-N/r},\qquad x\in{\mathbb R}^N,
  $$
  for some constant $C_1$;\vspace{3pt}
  \item $L^r\subset L^{r,\infty}$ and $L^r\not=L^{r,\infty}$ if $1<r<\infty$ and $L^{r,\infty}=L^r$ if $r\in\{1,\infty\}$;\vspace{3pt}
  \item Let $1<r<\infty$ and let $\{r_j\}_{j=1}^k\subset(1,\infty)$ be such that
  $$
  \frac{1}{r}=\frac{1}{r_1}+\cdots+\frac{1}{r_k}.
  $$
  Then there exists a constant $C_2$ such that
  $$
  \biggr\|\prod_{j=1}^kf_j\biggr\|_{r,\infty}\le C_2\prod_{j=1}^k\|f_j\|_{r_j,\infty}
  $$
  for $f_j\in L^{r_j,\infty}$ and $j=1,2,\dots,k$.
\end{itemize}
\vspace{8pt}

Now we state the main results of this paper but several explanations are in order. 
In most of the parabolic problem dealing with homogeneous equations, 
a crucial role is played by the fundamental solution. 
It happens that 
one can deduce several strong results as soon as one has an explicit form of the fundamental solution, 
allowing to get estimates for the function and its derivatives
(see for instance \cite{IKK01,IKK02, IKM, S}). 
In our problems, even if the coefficients are  rather simple, such an explicit form is unavailable. 
On the other hand,
bounds on the solution are known (see for instance \cite{CR, GN, GW}). 
In order to apply known bounds one has to impose additional properties on the weights under consideration.
More precisely, the weights have to belong to the $A_{1+\frac{2}{N}}$ class additionally to being $A_2$ 
and $w^{-N/2}$ has to satisfy a reverse doubling condition. 
We refer the reader to Section 2 for a discussion of these fact. 
In what follows, we put 
$$
p_*(\alpha):=1+\frac{2-\alpha}{N}\qquad\mbox{for}\quad \alpha \in \left \{ a,b \right \}.
$$
Furthermore, we assume that
\begin{gather}
\label{eq:1.4}
a\in(-1,1)\quad\mbox{if $N=1,2$},
\qquad
a\in(-1,2/N)\quad \mbox{if $N\ge3$}. \\ 
\label{eq:1.5}
b\in(-1,1)\quad\mbox{if $N=1$},
\qquad
b\in(-N,1)\quad \mbox{if $N\ge2$}.
\end{gather}

The first theorem is concerned with the nonexistence of global-in-time solutions of \eqref{eq:1.1}.

\begin{theorem}
\label{Theorem:1.1}
Assume \eqref{eq:1.4} and $1<p\le p_*(a)$.
Then problem \eqref{eq:1.1} with $w(x)=|x_1|^a$ has no nontrivial global-in-time solutions.
\end{theorem}
\begin{theorem}
\label{Theorem:1.2}
Assume \eqref{eq:1.5} and $1<p\le p_*(b)$.
Then problem \eqref{eq:1.1} with $w(x)=|x|^b$ has no nontrivial global-in-time solutions.
\end{theorem}

In second theorem we give a sufficient condition for 
the existence of nontrivial global-in-time solutions of \eqref{eq:1.1}.

\begin{theorem}
\label{Theorem:1.3}
Assume \eqref{eq:1.4} and $p>p_*(a)$.
Put
\begin{equation}
\label{eq:1.6}
r_*:=\frac{N}{2-a}(p-1)>1.
\end{equation}
Then the following holds:
\begin{itemize}
\item[\rm(i)]
There exists a positive constant $\delta$ such that,
for any $\varphi\in L^\infty\cap L^{r_*,\infty}$ with
\begin{equation}
\label{eq:1.7}
\|\varphi\|_{r_*,\infty}<\delta,
\end{equation}
a global-in-time solution~$u$ of \eqref{eq:1.1} with $w(x)=|x_1|^a$ exists and it satisfies
\begin{equation}
\label{eq:1.8}
\sup_{t>0}\,(1+t)^{\frac{N}{2-a}(\frac{1}{r_*}-\frac{1}{q})}\|u(t)\|_{q,\infty}<\infty,\qquad r_*\le q\le\infty.
\end{equation}
\item[\rm(ii)]
Let $1\le r\le r_*$.
Then there exists a positive constant $\delta$ such that,
for any $\varphi\in L^\infty\cap L^r$ with
\begin{equation}
\label{eq:1.9}
\|\varphi\|_r^{\frac{r}{r_*}}\|\varphi\|_\infty^{1-\frac{r}{r_*}}<\delta,
\end{equation}
a global-in-time solution~$u$ of \eqref{eq:1.1} with $w(x)=|x_1|^a$ exists and it satisfies
\begin{equation}
\label{eq:1.10}
\sup_{t>0}\,(1+t)^{\frac{N}{2-a}(\frac{1}{r}-\frac{1}{q})}\|u(t)\|_q<\infty,\qquad r\le q\le\infty.
\end{equation}
\end{itemize}
\end{theorem}
\begin{theorem}
\label{Theorem:1.4}
Assume \eqref{eq:1.5} and $p>p_*(b)$.
Put
$$
r_*:=\frac{N}{2-b}(p-1)>1.
$$
Then the following holds:
\begin{itemize}
\item[\rm(i)]
There exists a positive constant $\delta$ such that,
for any $\varphi\in L^\infty\cap L^{r_*,\infty}$ with
$$
\|\varphi\|_{r_*,\infty}<\delta,
$$
a global-in-time solution~$u$ of \eqref{eq:1.1} with $w(x)=|x|^b$ exists and it satisfies
\begin{equation}
\label{eq:1.11}
\sup_{t>0}\,(1+t)^{\frac{N}{2-b}(\frac{1}{r_*}-\frac{1}{q})}\|u(t)\|_{q,\infty}<\infty,\qquad r_*\le q\le\infty,
\end{equation}
\item[\rm(ii)]
Let $1\le r\le r_*$.
Then there exists a positive constant $\delta$ such that,
for any $\varphi\in L^\infty\cap L^r$ with
$$
\|\varphi\|_r^{\frac{r}{r_*}}\|\varphi\|_\infty^{1-\frac{r}{r_*}}<\delta.
$$
a global-in-time solution~$u$ of \eqref{eq:1.1} with $w(x)=|x|^b$ exists and it satisfies
$$
\sup_{t>0}\,(1+t)^{\frac{N}{2-b}(\frac{1}{r}-\frac{1}{q})}\|u(t)\|_q<\infty,\qquad r\le q\le\infty.
$$
\end{itemize}
\end{theorem}

As a direct consequence of Theorems~\ref{Theorem:1.3} and~\ref{Theorem:1.4}, we have:
\begin{corollary}
\label{Corollary:1.1}
Let $\alpha\in\{a,b\}$.
Assume $p>p_*(\alpha)$.
Then there exists a positive constant $\delta$ such that, if
  \begin{equation}
  \label{eq:1.12}
  |\varphi(x)|\le\frac{\delta}{1+|x|^{(2-\alpha)/(p-1)}},\qquad x\in{\mathbb R}^N,
  \end{equation}
  then a global-in-time solution $u$ of \eqref{eq:1.1} exists and it satisfies \eqref{eq:1.8} and \eqref{eq:1.11},
  respectively.
\end{corollary}
\begin{remark}
\label{Remark:1.1}
If we have suitable bounds for the derivatives of the fundamental solution,
then, applying the arguments in \cite{IKK01},
we can obtain the asymptotic behavior of solutions for \eqref{eq:1.1}.
However, unfortunately, it seems that they have been still left open. 
\end{remark}

The main technical difficulties arise in the case of $w(x)=|x_1|^a$, 
so we will give the proofs only in the first case, 
namely Theorems~\ref{Theorem:1.1} and ~\ref{Theorem:1.3}. 
However in the next section, we will check the restrictions on the range of exponents in the two cases. 

\section{Preliminaries}
A crucial tool in our arguments is based on the use of the fundamental solution of 
the operator $\partial_t -\text{div}(w(x)\nabla \cdot)$. 
As already mentioned due to the inhomogeneity of the operator, 
an explicit formula is not known but bounds are available (see below). 
In order to check these bounds one has to check  that the coefficient $w(x)$ is a $A_{1+\frac{2}{N}}$ weight 
in the sense of Muckenhoupt class 
and that the function $w^{-N/2}$ satisfies a doubling and reverse doubling condition 
of order $\mu$ with $\mu > 1/2$. 
Here we say that the function $w^{-N/2}$ satisfies doubling and a reverse doubling conditions if 
there exist positive constants $C_1$ and $C_2$ such that 
\[ 
	\int_{B_{sR}(x)}w(y)^{-\frac{N}{2}} dy \le C_1 s^{\mu N} \int_{B_{R}(x)}w(y)^{-\frac{N}{2}} dy
\] 
and 
\[ 
	\int_{B_{sR}(x)}w(y)^{-\frac{N}{2}} dy \ge C_2 s^{\mu N} \int_{B_{R}(x)}w(y)^{-\frac{N}{2}} dy
\] 
for all $x\in \mathbb{R}^N$, $s\ge 1$ and $R>0$, respectively. 

In the following we check these conditions for our model weights. 
Let $a\in(-1,1)$.
We put
\begin{equation}
\label{eq:2.1}
w(x)=|x_1|^a\ge0,
\end{equation}
and assume \eqref{eq:1.4}.
Then $w(x)$ is a $A_{1+\frac{2}{N}}$ weight in the sense of Muckenhoupt class.
Furthermore, the function $w^{-N/2}$ satisfies a doubling condition of order $1-a/2$
and reverse doubling condition of order $1-a/2$, 
thus $w^{-N/2}$ satisfies doubling and a reverse doubling condition of order $\mu$ with $\mu>1/2$ 
under the condition \eqref{eq:1.4}.

In the other case that $w(x)=|x|^b$, 
this is an $A_2$ weight as soon as $b \in (-N,N)$ and an $A_{1+\frac{2}{N}}$ weight as soon as $b \in (-N,2)$.
The doubling and reverse doubling conditions are checked for order $1-b/2$ 
and then one needs additionally that $b<1$ in order to check $w^{-N/2}$ satisfies a doubling 
and a reverse doubling condition of order $\mu$ with $\mu>1/2$.

 All in all these give the conditions on $a$ and $b$ described in the introduction. 
Under this situation,
the fundamental solution $\Gamma=\Gamma(x,y,t)$ has the following properties:
 \begin{itemize}
  	\item[(K1)] 
	$\displaystyle{\int_{{\mathbb R}^N}\Gamma(x,y,t)\, dx=\int_{{\mathbb R}^N}\Gamma(x,y,t)\, dy=1}$ 
	for $x,y\in{\mathbb R}^N$ and $t>0$;
	\item[(K2)]
	$\displaystyle{
	\Gamma(x,y,t)=\int_{{\mathbb R}^N}\Gamma(x,\xi,t-s)\,\Gamma(\xi,y,s)\,d\xi}
	$\quad
	for $x,y\in{\mathbb R}^N$ and $t>s>0$;
	\item[(K3)]
	There exist positive constants $c_*$ and $C_*$ such that 
	\begin{equation*}
	\begin{split}
	&
	c_*\left(\frac{1}{[h_x^{-1}(t)]^N}+\frac{1}{[h_y^{-1}(t)]^N}\right)
	e^{-c_*\left(\frac{h_x(|x-y|)}{t}\right)^\frac{1}{1-\alpha}}
	\le \Gamma(x,y,t)
	\\
	&\qquad\qquad\qquad\qquad
	\le C_*\left(\frac{1}{[h_x^{-1}(t)]^N}+\frac{1}{[h_y^{-1}(t)]^N}\right)
	e^{-C_*\left(\frac{h_x(|x-y|)}{t}\right)^\frac{1}{1-\alpha}}
	\end{split}
	\end{equation*}
		for $x,y\in{\mathbb R}^N$ and $t>0$, where $\alpha\in\{a,b\}$.
		Here
		\begin{equation}
		\label{eq:2.2}
		h_x(r)=\left(\int_{B_r(x)}w(y)^{-\frac{N}{2}}\,dy\right)^{\frac{2}{N}}
		\end{equation}
		and $h^{-1}_x$ denotes the inverse function of $h_x$.
  \end{itemize}
See \cite{GN}. (See also \cite{CR} and \cite{GW}.) 
From here, we focus on the case of $w(x)=|x_1|^a$.
By \eqref{eq:2.1} and \eqref{eq:2.2} we state a lemma on upper and lower estimates of $h_x(r)$.
In what follows,
by the letters $C$ and $C'$
we denote generic positive constants (independent of $x$ and $t$)
and they may have different values also within the same line. 
\begin{lemma}
\label{Lemma:2.1}
Let $a\in(-1,1)$. 
Then the following hold. 
\begin{itemize}
\item[\rm(i)]
For $a\in[0,1)$,
there exist positive constants $C$ and $C'$ such that
\begin{equation}
\label{eq:2.3}
h_x(r)\le Cr^{2-a}
\end{equation}
and
\begin{equation}
\label{eq:2.4}
h_x(r)\ge
C'
\left\{
\begin{array}{ll}
r^2|x_1|^{-a}
&
\mbox{if}\qquad 0<r\le |x_1|,
\\
\, &\\
r^{2-a}
&
\mbox{if}\qquad r\ge |x_1|,
\end{array}
\right.
\end{equation}
for all $x\in{\mathbb R}^N$ and $r>0$.
\item[\rm(ii)]
For $a\in(-1,0)$,
there exist positive constants $C$ and $C'$ such that
\begin{equation}
\label{eq:2.5}
h_x(r)\le
C
\left\{
\begin{array}{ll}
r^2|x_1|^{-a}
&
\mbox{if}\qquad 0<r\le |x_1|,
\\
\, &\\
r^{2-a}
&
\mbox{if}\qquad r\ge |x_1|,
\end{array}
\right.
\end{equation}
and
\begin{equation}
\label{eq:2.6}
h_x(r)\ge
C'r^{2-a}
\end{equation}
for all $x\in{\mathbb R}^N$ and $r>0$.
\end{itemize}
\end{lemma}
{\bf Proof.}
We first prove assertion~(i).
Since $w(y)^{-N/2}$ is monotonically decreasing function with respect to the distance from the origin,
by \eqref{eq:2.1} and \eqref{eq:2.2} we have
$$
h_x(r)=\left(\int_{B_r(x)}|y_1|^{-\frac{aN}{2}}\,dy\right)^{\frac{2}{N}}
\le \left(\int_{B_r(0)}|y_1|^{-\frac{aN}{2}}\,dy\right)^{\frac{2}{N}}
\le C r^{2-a}
$$
for all $x\in{\mathbb R}^N$ and $r>0$.
This implies \eqref{eq:2.3}.
On the other hand,
since $w(y)$ depends only on $y_1$ variable,
for any $x=(x_1,x')\in{\mathbb R}\times{\mathbb R}^{N-1}$,
we can choose a point $x_*=(x_1,0)$ such that
\begin{equation}
\label{eq:2.7}
\int_{B_r(x)}w(y)^{-\frac{N}{2}}\,dy=\int_{B_r(x_*)}w(y)^{-\frac{N}{2}}\,dy,\qquad r>0.
\end{equation}
Furthermore, for any $r>0$, we can take a cube $Q$ in ${\mathbb R}^N$ such that
\begin{equation}
\label{eq:2.8}
Q:=\left(x_1-\frac{r}{\sqrt{N}}, x_1+\frac{r}{\sqrt{N}}\right)\times 
\left(-\frac{r}{\sqrt{N}}, \frac{r}{\sqrt{N}}\right)^{N-1}\subset B_r(x_*).
\end{equation}
By \eqref{eq:2.1}, \eqref{eq:2.7} and \eqref{eq:2.8}, for any $x\in{\mathbb R}^N$ and $r>0$, we have
\begin{equation*}
\begin{split}
\int_{B_r(x)}w(y)^{-\frac{N}{2}}\,dy
&
=\int_{B_r(x_*)}|y_1|^{-\frac{aN}{2}}\,dy
\\
&
\ge \int_{Q}|y_1|^{-\frac{aN}{2}}\,dy
\\
&
\ge Cr^{N-1}\int_{x_1-\frac{r}{\sqrt{N}}}^{x_1+\frac{r}{\sqrt{N}}}|y_1|^{-\frac{aN}{2}}\, dy_1
\\
&
\ge Cr^N\min\left\{\left|x_1-\frac{r}{\sqrt{N}}\right|^{-\frac{aN}{2}}, \left|x_1+\frac{r}{\sqrt{N}}\right|^{-\frac{aN}{2}}\right\}
\\
&
\ge Cr^N(|x_1|+r)^{-\frac{aN}{2}}.
\end{split}
\end{equation*}
This together with \eqref{eq:2.2} yields \eqref{eq:2.4}.
Thus assertion~(i) holds.

Next we prove assertion~(ii).
Since $w(y)^{-N/2}$ is monotonically increasing function with respect to the distance from the origin,
by \eqref{eq:2.1} and \eqref{eq:2.2} we have
$$
h_x(r)=\left(\int_{B_r(x)}|y_1|^{-\frac{aN}{2}}\,dy\right)^{\frac{2}{N}}
\ge \left(\int_{B_r(0)}|y_1|^{-\frac{aN}{2}}\,dy\right)^{\frac{2}{N}}
\ge C r^{2-a}
$$
for all $x\in{\mathbb R}^N$ and $r>0$.
This implies \eqref{eq:2.6}.
On the other hand,
similarly to \eqref{eq:2.8},
for any $r>0$, we can take a cube $\tilde{Q}$ in ${\mathbb R}^N$ such that
\begin{equation}
\label{eq:2.9}
\tilde{Q}:=(x_1-r, x_1+r)\times(-r, r)^{N-1}\supset B_r(x_*).
\end{equation}
Then, since $a<0$,
by \eqref{eq:2.1}, \eqref{eq:2.7} and \eqref{eq:2.9}, for any $x\in{\mathbb R}^N$ and $r>0$, we have
\begin{equation*}
\begin{split}
\int_{B_r(x)}w(y)^{-\frac{N}{2}}\,dy
&
=\int_{B_r(x_*)}|y_1|^{-\frac{aN}{2}}\,dy
\\
&
\le \int_{\tilde{Q}}|y_1|^{-\frac{aN}{2}}\,dy
\\
&
\le Cr^{N-1}\int_{x_1-r}^{x_1+r}|y_1|^{-\frac{aN}{2}}\, dy_1
\\
&
\le Cr^N\max\{|x_1-r|^{-\frac{aN}{2}}, |x_1+r|^{-\frac{aN}{2}}\}
\\
&
\le Cr^N(|x_1|+r)^{-\frac{aN}{2}}.
\end{split}
\end{equation*}
This together with \eqref{eq:2.2} yields \eqref{eq:2.5}.
Thus assertion~(ii) holds, 
and Lemma~\ref{Lemma:2.1} follows.
\vspace{5pt}
$\Box$

For any $x\in{\mathbb R}^N$,
by \eqref{eq:2.1} and \eqref{eq:2.2} we can easily obtain that
\begin{equation}
\label{eq:2.10}
\frac{d}{dr}h_x(r)>0
\end{equation}
for all $r>0$.
Then, in the case $a\in[0,1)$ for $N=1,2$ or $a\in[0,2/N)$ for $N\ge3$, 
by \eqref{eq:2.3} and \eqref{eq:2.10}
we have
$$
h^{-1}_x(Cr^{2-a})\ge r
$$
for all $x\in{\mathbb R}^N$ and $r>0$,
and we see that
\begin{equation}
\label{eq:2.11}
h^{-1}_x(t)\ge Ct^{\frac{1}{2-a}}
\end{equation}
for all $x\in{\mathbb R}^N$ and $t>0$.
Similarly,
by \eqref{eq:2.4} we see that
\begin{equation*}
h^{-1}_x(t)\le
C
\left\{
\begin{array}{ll}
|x_1|^{\frac{a}{2}}t^{\frac{1}{2}}
&
\mbox{if}\qquad 0<t\le C|x_1|^{2-a},
\\
\, &\\
t^{\frac{1}{2-a}}
&
\mbox{if}\qquad t\ge C|x_1|^{2-a},
\end{array}
\right.
\end{equation*}
for all $x\in{\mathbb R}^N$ and $t>0$.
This together with Lemma~\ref{Lemma:2.1}, (K3) and \eqref{eq:2.11} implies that
\begin{equation}
\label{eq:2.12}
\begin{split}
&
C'\left(\min\{|x_1|^{-\frac{aN}{2}}, t^{-\frac{aN}{2}\frac{1}{2-a}}\}
+\min\{|y_1|^{-\frac{aN}{2}}, t^{-\frac{aN}{2}\frac{1}{2-a}}\}\right)
t^{-\frac{N}{2}}e^{-C'\left(\frac{|x-y|^{2-a}}{t}\right)^\frac{1}{1-a}}
\\
&\qquad\qquad\qquad\qquad\qquad\qquad\qquad\qquad\qquad\quad
\le \Gamma(x,y,t)
\le Ct^{-\frac{N}{2-a}}e^{-C\left(\frac{|x-y|^{2-a}}{t}\right)^\frac{1}{1-a}}
\end{split}
\end{equation}
for $x,y\in{\mathbb R}^N$ and $t>0$.
Similarly to \eqref{eq:2.12},
for the case $a\in(-1,0]$,
we see that
\begin{equation}
\label{eq:2.13}
\begin{split}
&
C't^{-\frac{N}{2-a}}e^{-C'\left(\frac{|x-y|^{2-a}}{t}\right)^\frac{1}{1-a}}
\le \Gamma(x,y,t)
\\
&\quad
\le C\left(\max\{|x_1|^{-\frac{aN}{2}}, t^{-\frac{aN}{2}\frac{1}{2-a}}\}
+\max\{|y_1|^{-\frac{aN}{2}}, t^{-\frac{aN}{2}\frac{1}{2-a}}\}\right)
t^{-\frac{N}{2}}e^{-C\left(\frac{|x-y|^{2-a}}{t}\right)^\frac{1}{1-a}}
\end{split}
\end{equation}
for $x,y\in{\mathbb R}^N$ and $t>0$.
Then, by (K1), \eqref{eq:2.12} and \eqref{eq:2.13},
for any $1\le r\le \infty$, we have
\begin{equation}
\label{eq:2.14}
\|\Gamma(y,t)\|_r
\le Ct^{-\frac{N}{2-a}(1-\frac{1}{r})}
\end{equation}
for all $y\in{\mathbb R}^N$ and $t>0$.
For any $\varphi\in L^\infty$,
we put
\begin{equation}
\label{eq:2.15}
[S(t)\varphi](x):=\int_{{\mathbb R}^N}\Gamma(x,y,t)\varphi(y)\,dy
\end{equation}
for all $x\in{\mathbb R}^N$ and $t>0$.
This together with the Young inequality
(see e.g. \cite[Section~4, Chapter~IX]{RS} and \cite[Theorem~2.10.1]{Z})
with \eqref{eq:2.14} implies that
\begin{itemize}
\item[\rm (G1)]
It holds that
$$
\|S(t)\varphi\|_r\le Ct^{-\frac{N}{2-a}(\frac{1}{q}-\frac{1}{r})}\|\varphi\|_q,\qquad t>0,
$$
for any $\varphi\in L^q$ and $1\le q\le r\le \infty$;
\item[\rm (G2)]
It holds that
$$
\|S(t)\varphi\|_{r,\infty}\le Ct^{-\frac{N}{2-a}(\frac{1}{q}-\frac{1}{r})}\|\varphi\|_{q,\infty},\qquad t>0,
$$
for any $\varphi\in L^{q,\infty}$ and $1\le q\le r\le \infty$.
\end{itemize}

Furthermore, by \eqref{eq:2.12} and \eqref{eq:2.13} we have the following lemmas.
\begin{lemma}
\label{Lemma:2.2}
Assume \eqref{eq:1.4}.
Let $\varphi\in L^1$ be a non-trivial measurable function such that 
$\varphi\ge 0$ in ${\mathbb R}^N$ 
and $\operatorname{supp}\varphi$ is compact.
Then there exist positive constants $C$ and $T$ such that 
$$
[S(t)\varphi](x)\ge Ct^{-\frac{N}{2-a}}\int_{{\mathbb R}^N}\varphi(y)\, dy
$$
for $|x|\le t^{\frac{1}{2-a}}$ and $t\ge T$.
\end{lemma}
{\bf Proof.}
By \eqref{eq:2.12} and \eqref{eq:2.13} we can find positive constants $C$ and $T$ such that 
$$
\Gamma(x,y,t)\ge C t^{-\frac{N}{2-a}}
$$
for $|x|\le t^{\frac{1}{2-a}}$, $y\in\mbox{supp}\,\varphi$ and $t\ge T$.
Then Lemma~\ref{Lemma:2.2} follows from \eqref{eq:2.15}.
$\Box$
\vspace{3pt}

\section{Proof of Theorem~\ref{Theorem:1.1}}
In this section we prove Theorem~\ref{Theorem:1.1}, 
which means that problem~\eqref{eq:1.1} has no nontrivial global solutions in the case $1<p\le p_*(a)$. 
The proof of Theorem~\ref{Theorem:1.1} is based on the arguments
\cite[Theorem~5]{W1}, and \cite[Theorem~1]{W2}
(see also \cite[Theorem~1.1]{FIK2}).
\vspace{5pt}

We first prove the following lemma. 
\begin{lemma}
\label{Lemma:3.1} 
Let $u$ be a solution of \eqref{eq:1.1} in ${{\mathbb R}^N}\times(0,T)$ with $0<T\le\infty$.
Then there exists a constant $C_*$, independent of $\varphi$ and $T$, such that
\begin{equation}
\label{eq:3.1}
t^{\frac{1}{p-1}}\|S(t)\varphi\|_\infty\le C_*
\end{equation}
for any  $t\in[0,T)$.
\end{lemma}
\noindent{\bf Proof.}
This lemma follows from the proof of \cite[Theorem~5]{W1}.
For completeness of this paper, we will add the proof of it.

Since it follows from \eqref{eq:2.12} and \eqref{eq:2.13} that the fundamental solution $\Gamma$ is positive
for $x,y\in{\mathbb R}^N$ and $t>0$,
by \eqref{eq:1.2} and \eqref{eq:2.15} we have
\begin{equation}
\label{eq:3.2}
[S(t)\varphi](x)\le u(x,t)<\infty 
\end{equation}
for almost all $x\in{\mathbb R}^N$ all $t\in(0,T)$. 
This together with \eqref{eq:1.2} again implies 
\begin{equation}
\label{eq:3.3}
u(x,t)\ge
\int_0^t\int_{{\mathbb R}^N}\Gamma(x,y,t-s)\left([S(s)\varphi](y)\right)^p\,dy\,ds
\end{equation}
for almost all $x\in{\mathbb R}^N$ and all $t\in(0,T)$.  
Then, applying  the Jensen inequality with the aid of $(K1)$ and $(K2)$ to \eqref{eq:3.3}, 
we obtain 
\begin{equation}
\label{eq:3.4}
u(x,t)
\ge 
\int_0^t
\biggr(\int_{{\mathbb R}^N}\Gamma(x,y,t-s)[S(s)\varphi](y)dy\biggr)^p\,ds
=t([S(t)\varphi](x))^p
\end{equation}
for almost all $x\in{\mathbb R}^N$ and all $t\in(0,T)$.  
We repeat the above argument with \eqref{eq:3.2} replaced by \eqref{eq:3.4}, 
and have 
\begin{equation*}
\begin{split}
u(x,t)
& \ge
\int_0^t\int_{{\mathbb R}^N}\Gamma(x,y,t-s)\bigg(s([S(s)\varphi](y))^p\bigg)^p\,dy\,ds
\\
&
\ge
\int_0^ts^p\biggr(\int_{{\mathbb R}^N}\Gamma(x,y,t-s)[S(s)\varphi](y)dy\biggr)^{p^2}\,ds
=\frac{1}{p+1}t^{p+1}([S(t)\varphi](x))^{p^2}
\end{split}
\end{equation*}
for almost all $x\in{\mathbb R}^N$ and all $t\in(0,T)$. 
Repeating the above argument,
for any $k=2,3,\dots$,
it holds that
\begin{equation}
\label{eq:3.5}
u(x,t)\ge A_kt^{\frac{p^k-1}{p-1}}\left([S(t)\varphi](x)\right)^{p^k}
\end{equation}
for almost all $x\in{\mathbb R}^N$ and all $t\in(0,T)$, 
where
\begin{equation*}
\begin{split}
 A_k:&=
 \left(\frac{1}{p+1}\right)^{p^{k-2}}\left(\frac{1}{(p+1)p+1}\right)^{p^{k-3}}
\cdots\left(\frac{1}{(1+p+\cdots+p^{k-2})p+1}\right)\\
&
=\prod_{j=1}^{k-1}\left(\frac{p-1}{p^{j+1}-1}\right)^{p^{k-j-1}}. 
\end{split}
\end{equation*}
Therefore, by \eqref{eq:3.5} we have 
\begin{equation}
\label{eq:3.6}
t^{\frac{1}{p-1}(1-\frac{1}{p^k})}
[S(t)\varphi](x)
\le
u(x,t)^{p^{-k}}\left(\prod_{j=1}^{k-1}\left(\frac{p-1}{p^{j+1}-1}\right)^{p^{k-j-1}}\right)^{-p^{-k}}
\end{equation}
for almost all $x\in{\mathbb R}^N$ and all $t\in(0,T)$. 
On the other hand, 
we have 
\begin{equation}
\label{eq:3.7}
\begin{split}
\log\left(\prod_{j=1}^{\infty}\left(\frac{p^{j+1}-1}{p-1}\right)^{p^{-j-1}}\right)
&=
\sum_{j=1}^\infty p^{-j-1}\log\left(\frac{p^{j+1}-1}{p-1}\right)
\\
&\le
\sum_{j=1}^\infty p^{-j-1}\log\left((j+1)p^j\right)<\infty.
\end{split}
\end{equation}
Then, by \eqref{eq:3.6} and \eqref{eq:3.7} 
we can find a constant $C_*$, independent of $T$ and the initial function $\varphi$, 
such that 
$$
t^{\frac{1}{p-1}}[S(t)\varphi](x)\le C_*<\infty
$$
for almost all $x\in{\mathbb R}^N$ and all $t\in(0,T)$. 
This implies \eqref{eq:3.1}, and Lemma~\ref{Lemma:3.1} follows. 
$\Box$\vspace{5pt}

We prove Theorem~\ref{Theorem:1.1} by using Lemma~\ref{Lemma:3.1}.
\vspace{5pt}
\newline
\noindent
{\bf Proof of Theorem~\ref{Theorem:1.1}.} 
The proof is by contradiction. 
Let $u$ be a global-in-time solution of \eqref{eq:1.1}. 
Since $u(\cdot,1)$ is a positive measurable function in ${\mathbb R}^N$, 
we can find a non-trivial measurable function $\varphi_1\in L^1$ such that 
$\operatorname{supp}\varphi_1$ is compact and
$0\le\varphi_1(x)\le u(x,1)$ for almost all $x\in{\mathbb R}^N$.
Then it follows from Lemma~\ref{Lemma:2.2} that 
\begin{equation}
\label{eq:3.8}
[S(t)\varphi_1](x)
\ge CMt^{-\frac{N}{2-a}},\qquad 
M:=\int_{{\mathbb R}^N}\varphi_1(x)\,dx,
\end{equation}
for all $|x| \leq t^{\frac{1}{2-a}}$ and $t\ge T$,
where $T$ is a positive constant given in Lemma~\ref{Lemma:2.2}. 
Furthermore, by \eqref{eq:1.2}, \eqref{eq:2.15} and $(K2)$ we see that 
\begin{equation}
\label{eq:3.9}
u(x,t+1)\ge[S(t)u(1)](x)\ge[S(t)\varphi_1](x)
\end{equation}
for almost all $x\in{\mathbb R}^N$ and all $t>0$.

We first consider the case $1<p<p_*(a)$.
By \eqref{eq:3.8} and \eqref{eq:3.9} we have
\begin{equation}
\label{eq:3.10}
[S(t)u(1)](x)\ge CMt^{-\frac{N}{2-a}}
\end{equation}
for all $|x|\leq t^{\frac{1}{2-a}}$ and $t\ge T$.
It follows from $1<p<p_*(a)$ with \eqref{eq:3.10} that
$$
t^{\frac{1}{p-1}}\|S(t)u(1)\|_\infty\to\infty\quad\mbox{as}\quad t\to\infty,
$$
which contradicts \eqref{eq:3.1}. 
This means that  problem~\eqref{eq:1.1} possesses no global-in-time positive solutions.

Next we consider the case $p=p_*(a)$.
By \eqref{eq:1.2}, \eqref{eq:3.8} and \eqref{eq:3.9} we have
\begin{equation}
\label{eq:3.11}
\begin{split}
 & \int_{{\mathbb R}^N}u(x,t+1)\,dx
\ge\int_{{\mathbb R}^N}\int_2^{t+1}\int_{{\mathbb R}^N}\Gamma(x,y,t+1-s)u(y,s)^p\,dy\,ds\,dx\\
 & 
=\int_2^{t+1}\int_{{\mathbb R}^N}u(y,s)^p\,dy\,ds
\ge\int_1^t\int_{\{s^{1/(2-a)}\le |y|\le 2s^{1/(2-a)}\}}u(y,s+1)^p\,dy\,ds\\
 & 
\ge C\int_1^t s^{-\frac{N}{2-a}(p-1)}\,ds=C\log t,\qquad t>1.
\end{split}
\end{equation}
Let $m$ be a sufficiently large positive constant. 
By \eqref{eq:3.11} we can find $T>0$ and 
a non-trivial measurable function $\varphi_2\in L^1$ such that 
$\operatorname{supp}\varphi_2$ is compact and
$$
0\le\varphi_2(x)\le u(x,T)\quad\mbox{for almost all $x\in{\mathbb R}^N$},
\qquad
\int_{{\mathbb R}^N}\varphi_2(x)\,dx\ge m. 
$$
Similarly to \eqref{eq:3.8} and \eqref{eq:3.9}, we have 
$$
u(x,t+T)\ge [S(t)\varphi_2](x)
\ge Cmt^{-\frac{N}{2-a}}
$$
for almost all $|x|\le t^{\frac{1}{2-a}}$ and all $t\ge 1$.
This implies that
\begin{equation}
\label{eq:3.12}
t^{\frac{N}{2-a}}\|S(t)\varphi_2\|_\infty\ge Cm,\qquad t>1.
\end{equation} 
Let $v$ be a solution of \eqref{eq:1.1} with initial data $\varphi_2$.
Then, since $u$ is a global-in-time solution of \eqref{eq:1.1}, 
$v$ is also a global-in-time solution of \eqref{eq:1.1}.
Therefore we can apply Lemma~\ref{Lemma:3.1} to the solution $v$,
and obtain \eqref{eq:3.1} replacing $\varphi$ with $\varphi_2$.
By the arbitrariness of $m$, 
this contradicts \eqref{eq:3.12}
and we see that problem~\eqref{eq:1.1} possesses no global-in-time positive solutions for the case $p=p_*(a)$.
Therefore the proof of Theorem~\ref{Theorem:1.1} is complete.
$\Box$

\section{Proof of Theorem~\ref{Theorem:1.3}}
In this section we prove Theorem~\ref{Theorem:1.3}.
We first prove the uniqueness of solutions of \eqref{eq:1.1}. 
(See also \cite[Lemma~3.1]{FIK}.)

\begin{lemma}
\label{Lemma:4.1}
Let $i=1,2,\tau>0$, and $u_i$ be a solution of \eqref{eq:1.1} in ${{\mathbb R}^N}\times(0,\tau)$ with
$\varphi=\varphi_i\in L^\infty$.
Then,
for any $\sigma\in(0,\tau)$,
there exists a constant $C$ such that
$$
\sup_{0<t\le\sigma}\|u_1(t)-u_2(t)\|_\infty\le C\|\varphi_1-\varphi_2\|_\infty.
$$
Here the constant $C$ depends on $\|u_1\|_{L^\infty(0,\sigma:L^\infty)}$
and $\|u_2\|_{L^\infty(0,\sigma:L^\infty)}$.
\end{lemma}
{\bf Proof.}
Let $\sigma\in(0,\tau)$.
Put $v=u_1-u_2$.
Then we have
$$
\|v\|_{L^\infty(0,\sigma:L^\infty)}\le\|u_1\|_{L^\infty(0,\sigma:L^\infty)}
+\|u_2\|_{L^\infty(0,\sigma:L^\infty)}<\infty.
$$
This together with \eqref{eq:1.2} and (K2) yields
\begin{equation*}
\begin{split}
|v(x,\tilde{t})|
&\le
\|v(t)\|_\infty+\int_t^{\tilde{t}}\int_{{\mathbb R}^N}\Gamma(x,y,\tilde{t}-s)|u_1(y,s)^p-u_2(y,s)^p|\,dy\, ds\\
&\le
\|v(t)\|_\infty+C_1\int_t^{\tilde{t}}\int_{{\mathbb R}^N}\Gamma(x,y,\tilde{t}-s)|v(y,s)|\,dy\, ds\\
&\le
\|v(t)\|_\infty+C_1\sup_{t<\tau\le \tilde{t}}\|v(\tau)\|_\infty(\tilde{t}-t)
\end{split}
\end{equation*}
for almost all $x\in {{\mathbb R}^N}$ and all $0\le t<\tilde{t}\le\sigma$,
where $C_1$ is a positive constant.
This implies that
$$
\sup_{t<\tau\le \tilde{t}}\|v(t)\|_\infty\le\|v(t)\|_\infty+C_1\sup_{t<\tau\le \tilde{t}}\|v(\tau)\|_\infty(\tilde{t}-t)
$$
for all $0\le t<\tilde{t}\le\sigma$.

Let $\varepsilon$ be a sufficiently small positive constant such that $C_1\varepsilon\le1/2$ and
$\varepsilon<\sigma$.
Then, by \eqref{eq:3.2} we have
$$
\sup_{t<\tau\le t+\varepsilon}\|v(\tau)\|_\infty\le2\|v(t)\|_\infty
$$
for all $t\in[0,\sigma-\varepsilon]$.
Therefore there exists a constant $C_2$ such that
$$
\sup_{0<\tau\le\sigma}\|v(t)\|_\infty\le C_2\|v(0)\|_\infty,
$$
and we have inequality \eqref{eq:3.1}.
Thus the proof of Lemma \ref{Lemma:3.1} is complete.
$\Box$
\vspace{5pt}

Next we prove local existence of solutions of \eqref{eq:1.2}.
For any nonnegative function $\varphi\in L^\infty$,
we define $\{u_n\}$ inductively by
\begin{equation}
\label{eq:4.1}
\begin{array}{l}
\displaystyle
\qquad
u_1(x,t):=\int_{{\mathbb R}^N}\Gamma(x,y,t)\varphi(y)\,dy,
\vspace{5pt}\\
\displaystyle
\qquad
u_{n+1}(x,t):=u_1(x,t)+\int_0^t\int_{{\mathbb R}^N}\Gamma(x,y,t-s)u_n(y,s)^p\,dy\,ds,
\qquad n=1,2,\cdots,
\end{array}
\end{equation}
for almost all $x\in{\mathbb R}^N$ and all $t>0$.
Then we can prove inductively that
\begin{equation}
\label{eq:4.2}
0\le u_{n-1}(x,t)\le u_n(x,t)
\end{equation}
for almost all $x\in{\mathbb R}^N$ and all $t>0$.
This means that the limit function
\begin{equation}
\label{eq:4.3}
u_*(x,t):=\lim_{n\to\infty}u_n(x,t)\in[0,\infty]
\end{equation}
can be defined for almost all $x\in{\mathbb R}^N$ and all $t>0$.
Furthermore,
by properties (G1) and (G2) we can find a constant $c_*$ such that 
\begin{equation}
\label{eq:4.4}
\begin{split}
 & \sup_{0<t<\infty}\|u_1(t)\|_\infty\le c_*\|\varphi\|_\infty,\\
 & \sup_{0<t<\infty}t^{\frac{N}{2-a}(\frac{1}{r_*}-\frac{1}{q})}\|u_1(t)\|_{q,\infty}
 \le c_*\|\varphi\|_{r_*,\infty},
 \end{split}
\end{equation}
for any $q\in[r_*,\infty]$ if $\varphi\in L^{r_*,\infty}$.
The we have the following.
(See also \cite[Lemma~3.2]{FIK} and \cite[Lemma~3.1]{IKK01}.)
\begin{lemma}
\label{Lemma:4.2}
Let $\varphi\in L^\infty$.
Then there exists a solution~$u$  of \eqref{eq:1.1} in ${\mathbb R}^N\times(0,T)$ for some $T>0$ 
such that 
\begin{equation}
\label{eq:4.5}
 \sup_{0<t<T}\|u(t)\|_\infty\le 2c_*\|\varphi\|_\infty.
\end{equation}
Here $c_*$ is the constant given in \eqref{eq:4.4}. 
\end{lemma}
{\bf Proof.}
Let $T$ be a sufficiently small positive constant to be chosen later. 
By induction we prove that 
\begin{equation}
\label{eq:4.6}
 \sup_{0<t<T}\|u_n(t)\|_\infty\le 2c_*\|\varphi\|_\infty
\end{equation}
for all $n=1,2,\dots$.
By \eqref{eq:4.4} we have \eqref{eq:4.6} for $n=1$. 
Assume that \eqref{eq:4.6} holds for some $n=n_*\in\{1,2,\dots\}$, that is, 
$$
\sup_{0<t<T}\|u_{n_*}(t)\|_\infty\le 2c_*\|\varphi\|_\infty.
$$
Then, by \eqref{eq:4.1} and (G1) we have
\begin{equation}
\label{eq:4.7}
\begin{split}
\|u_{n_*+1}(t)\|_\infty
& \le
\|u_1(t)\|_\infty+\int_0^t \|S(t-s)u_{n_*}(s)^p\|_\infty\, ds
\\
& \le 
c_*\|\varphi\|_\infty+C_1\int_0^t\|u_{n_*}(s)\|_\infty^p\,ds
\\
& \le 
c_*\|\varphi\|_\infty
+C_1T(2c_*\|\varphi\|_\infty)^p
\end{split}
\end{equation}
for all $t\in(0,T)$, where $C_1$ is a constant independent of $n_*$ and $T$. 
Let $T$ be a sufficiently small constant such that 
$$
C_1T2^p(c_*\|\varphi\|_\infty)^{p-1}\le 1.
$$
Then, by \eqref{eq:4.7} we have \eqref{eq:4.6} for $n=n_*+1$. 
Therefore \eqref{eq:4.6} holds for all $n=1,2,\dots$. 
By \eqref{eq:4.2}, \eqref{eq:4.3} and \eqref{eq:4.6} we see that
the limit function $u_*$ satisfies \eqref{eq:1.2} and
$$
 \sup_{0<t<T}\|u_*(t)\|_\infty\le 2c_*\|\varphi\|_\infty.
$$
This together with Lemma~\ref{Lemma:4.1} implies that 
the function $u=u_*$ is a solution of \eqref{eq:1.1} in ${\mathbb R}^N\times(0,T)$. 
Thus Lemma~\ref{Lemma:4.2} follows.
$\Box$
\vspace{5pt}

Now we are ready to prove Theorem~\ref{Theorem:1.3}.
\vspace{5pt}

\noindent
{\bf Proof of the assertion~(i) of Theorem~\ref{Theorem:1.3}.}
Let $\delta$ be a sufficiently small positive constant. 
Assume \eqref{eq:1.7}. 
By induction we prove
\begin{equation}
\label{eq:4.8}
\begin{split}
&
\|u_n(t)\|_{r_*,\infty}\le 2c_*\delta,
\\
&
\|u_n(t)\|_\infty\le 2c_*\delta t^{-\frac{N}{(2-a)r_*}},
\end{split}
\end{equation}
for all $t>0$. 
By \eqref{eq:4.4} we have \eqref{eq:4.8} for $n=1$. 
Assume that \eqref{eq:4.8} holds for some $n=n_*\in\{1,2,\dots\}$, that is, 
\begin{equation}
\label{eq:4.9}
\begin{split}
&
\|u_{n_*}(t)\|_{r_*,\infty}\le 2c_*\delta,
\\
&
\|u_{n_*}(t)\|_\infty\le 2c_*\delta t^{-\frac{N}{(2-a)r_*}},
\end{split}
\end{equation}
for all $t>0$. 
These imply that 
\begin{equation}
\label{eq:4.10}
\|u_{n_*}(t)\|_{q,\infty}
\le \|u_{n_*}(t)\|_{r_*, \infty}^{\frac{r_*}{q}}\|u_{n_*}(t)\|_{\infty}^{1-\frac{r_*}{q}}
\le 2c_*\delta t^{-\frac{N}{2-a}(\frac{1}{r_*}-\frac{1}{q})}
\end{equation}
for all $t>0$ and $r_*<q<\infty$.
Since $r_*=N(p-1)/(2-a)$,
by \eqref{eq:4.9} we have
\begin{equation}
\label{eq:4.11}
\|u_{n_*}(t)^p\|_\infty
=\|u_{n_*}(t)\|_\infty^p
\le \left( 2c_*\delta t^{-\frac{N}{(2-a)r_*}}\right)^p
=(2c_*\delta)^pt^{-\frac{N}{(2-a)r_*}-1}
\end{equation}
for all $t>0$.
Similarly,
for any $\eta>1$ with $\eta\le r_*<\eta p$, 
by \eqref{eq:4.10}
we obtain
\begin{equation}
\label{eq:4.12}
\|u_{n_*}(t)^p\|_{\eta,\infty} 
=\|u_{n_*}(t)\|_{\eta p,\infty}^p
\le \left(2c_*\delta t^{-\frac{N}{2-a}(\frac{1}{r_*}-\frac{1}{\eta p})}\right)^p
\le C\delta^pt^{\frac{N}{(2-a)\eta}-\frac{N}{(2-a)r_*}-1}
\end{equation}
for all $t>0$.
Therefore, by (G1), (G2), \eqref{eq:4.11} and \eqref{eq:4.12}
we have  
\begin{equation}
\label{eq:4.13}
\begin{split}
 \left\|\int_{t/2}^t  S(t-s)u_{n_*}(s)^p\,ds\right\|_\infty
 &\le\int_{t/2}^t\|S(t-s)u_{n_*}(s)^p\|_\infty\,ds
 \\
 &
 \le \int_{t/2}^t\|u_{n_*}(s)^p\|_\infty\,ds\\
 &
 \le C\delta^p\int_{t/2}^ts^{-\frac{N}{(2-a)r_*}-1}\,ds
 \le C\delta^pt^{-\frac{N}{(2-a)r_*}}
\end{split}
\end{equation}
and 
\begin{equation}
\label{eq:4.14}
\begin{split}
 \left\|\int_{t/2}^t  S(t-s)u_{n_*}(s)^p\,ds\right\|_{r_*,\infty}
 &\le\int_{t/2}^t\|S(t-s)u_{n_*}(s)^p\|_{r_*,\infty}\,ds
 \\
 &
 \le \int_{t/2}^t\|u_{n_*}(s)^p\|_{r_*,\infty}\,ds
 \le C\delta^p\int_{t/2}^ts^{-1}\,ds
 \le C\delta^p
\end{split}
\end{equation}
for all $t>0$. 
On the other hand, 
by (G2), \eqref{eq:4.11} and \eqref{eq:4.12} with $\eta<r_*$ 
we have 
\begin{equation}
\label{eq:4.15}
\begin{split}
 &\left\|\int_0^{t/2} S(t-s)u_{n_*}(s)^p\,ds\right\|_\infty
 \\
 & \le\int_0^{t/2}\|S(t-s)u_{n_*}(s)^p\|_\infty\,ds
\le C\int_0^{t/2}(t-s)^{-\frac{N}{(2-a)\eta}}\|u_{n_*}(s)^p\|_{\eta,\infty}ds
\\
 &
 \le C \delta^pt^{-\frac{N}{(2-a)\eta}}\int_0^{t/2}s^{\frac{N}{(2-a)\eta}-\frac{N}{(2-a)r_*}-1}\,ds
 \le C\delta^pt^{-\frac{N}{(2-a)r_*}}
 \end{split}
\end{equation}
and 
\begin{equation}
\label{eq:4.16}
\begin{split}
 & \left\|\int_0^{t/2}S(t-s)u_{n_*}(s)^p\,ds\right\|_{r_*,\infty}
 \\
 & \le\int_0^{t/2}\|S(t-s)u_{n_*}(s)^p\|_{r_*,\infty}\,ds
\le C\int_0^{t/2}(t-s)^{-\frac{N}{2-a}(\frac{1}{\eta}-\frac{1}{r_*})}\|u_{n_*}(s)^p\|_{\eta,\infty}\,ds
\\
 &
 \le C \delta^pt^{-\frac{N}{2-a}(\frac{1}{\eta}-\frac{1}{r_*})}\int_0^{t/2}s^{\frac{N}{(2-a)\eta}-\frac{N}{(2-a)r_*}-1}\,ds
\le C \delta^p
 \end{split}
\end{equation}
for all $t>0$. 
Then, taking a sufficiently small $\delta$ if necessary, 
by \eqref{eq:4.4} and \eqref{eq:4.13}, \eqref{eq:4.14}, \eqref{eq:4.15} and \eqref{eq:4.16} 
we see that  
$$
\left.
\begin{array}{l}
t^{\frac{N}{(2-a)r_*}}\|u_{n_*+1}(t)\|_\infty
\vspace{5pt}\\
\quad\,\,\,\,\,
\|u_{n_*+1}(t)\|_{r_*,\infty}
\end{array}
\right\}
\le c_*\delta+C_1\delta^p\le 2c_*\delta
$$
for all $t>0$, where $C_1$ is a constant independent of $n_*$ and $\delta$. 
Hence we obtain \eqref{eq:4.8} for $n=n_*+1$. 
Thus \eqref{eq:4.8} holds for all $n=1,2,\dots$. 
Therefore, 
applying a similar argument as in the proof of Lemma~\ref{Lemma:4.2}, 
by \eqref{eq:4.8} we see that 
there exists a global-in-time solution $u$ of \eqref{eq:1.1} such that 
$$
\|u(t)\|_{r_*,\infty}\le 2c_*\delta,\qquad
\|u(t)\|_\infty\le 2c_*\delta t^{-\frac{N}{(2-a)r_*}},
$$
for all $t>0$.
This together with \eqref{eq:4.5} implies that
$$
\|u(t)\|_\infty\le C(1+t)^{-\frac{N}{(2-a)r_*}}
$$
 for all $t>0$.
 Furthermore, we apply an interpolation theorem to obtain
$$ 
\|u(t)\|_{q,\infty}\le C(1+t)^{-\frac{N}{2-a}(\frac{1}{r_*}-\frac{1}{q})},\qquad r_*\le q\le\infty,
$$
for all $t>0$. 
Thus we have \eqref{eq:1.8}, and the proof of the assertion of Theorem~\ref{Theorem:1.3} is complete. 
$\Box$
\vspace{5pt}

\noindent{\bf Proof of the assertion~(ii) of Theorem~\ref{Theorem:1.3}.}
Let $\delta$ be a sufficiently small constant and assume \eqref{eq:1.9}. 
Then, by the assertion~(i) of Theorem~\ref{Theorem:1.3} we see that 
there exists a global-in-time solution $u$ of \eqref{eq:1.1} satisfying \eqref{eq:1.8}. 

We prove the existence of a global-in-time solution of \eqref{eq:1.1} satisfying \eqref{eq:1.10}. 
For $r=r_*$, it follows from a similar argument as in the proof of the assertion~(i) of Theorem~\ref{Theorem:1.3}. 
So we assume $1\le r<r_*$. 
Put 
$$
\varphi_\lambda(x):=\lambda^\alpha\varphi(\lambda x),
\qquad
u_{n,\lambda}(x,t):=\lambda^\alpha u_n(\lambda x,\lambda^{2-a} t),
$$
where $\alpha=N/r_*$ and $\lambda$ is a positive constant such that 
$$
\|\varphi_\lambda\|_r=\|\varphi_\lambda\|_\infty. 
$$
Since 
$$
\|\varphi_\lambda\|_r^{\frac{r}{r_*}}\|\varphi_\lambda\|_\infty^{1-\frac{r}{r_*}}=\|\varphi\|_r^{\frac{r}{r_*}}\|\varphi\|_\infty^{1-\frac{r}{r_*}},
$$ 
it follows from \eqref{eq:1.9} that 
\begin{equation}
\label{eq:4.17}
\|\varphi_\lambda\|_r=\|\varphi_\lambda\|_\infty<\delta. 
\end{equation}
Furthermore, $u_{n,\lambda}$ satisfies 
\begin{equation}
\label{eq:4.18}
u_{n,\lambda}(t)
=S(t-\tau)u_{n,\lambda}(\tau)+\int_{\tau}^t S(t-s)u_{n-1,\lambda}(s)^p\,ds,
\end{equation}
for all $t>\tau\ge 0$.
On the other hand, by (G1), \eqref{eq:4.4} and \eqref{eq:4.17}
we can find a constant $C_*$ independent of $\delta$ such that 
\begin{equation}
\label{eq:4.19}
\|S(t)\varphi_\lambda\|_q\le C_*\delta(1+t)^{-\frac{N}{2-a}(\frac{1}{r}-\frac{1}{q})},\qquad t>0,
\end{equation}
for any $q\in[r,\infty]$.

By induction we prove that 
\begin{equation}
\label{eq:4.20}
\|u_{n,\lambda}(t)\|_q\le 2C_*\delta,
\quad 0<t\le 2, 
\end{equation}
for any $q\in[r,\infty]$ and $n=1,2,\dots$. 
By \eqref{eq:4.19} we have \eqref{eq:4.20} for $n=1$. 
Assume that \eqref{eq:4.20} holds for some $n=n_*$, that is, 
\begin{equation}
\label{eq:4.21}
\|u_{n_*,\lambda}(t)\|_q\le 2C_*\delta,
\quad 0<t\le 2, 
\end{equation}
for any $q\in[r,\infty]$.
Then, by \eqref{eq:4.21}, for any $q\in[r,\infty]$, we have 
\begin{equation}
\label{eq:4.22}
\|u_{n_*,\lambda}(t)^p\|_q
= \|u_{n_*,\lambda}(t)\|_{pq}^p
\le (2C_*\delta)^p
\end{equation}
for all $0<t\le 2$.
Taking a sufficiently small $\delta$ if necessary, 
by (G1), \eqref{eq:4.18}, \eqref{eq:4.19} and \eqref{eq:4.22} 
we obtain  
\begin{equation}
\label{eq:4.23}
\begin{split}
\|u_{n_*+1,\lambda}(t)\|_q 
& 
\le \|S(t)\varphi_\lambda\|_q+\int_0^t\|S(t-s)u_{n_*,\lambda}(s)^p\|_q\,ds
\\
& 
\le \|S(t)\varphi_\lambda\|_q+C_1\int_0^t\|u_{n_*,\lambda}(s)^p\|_q\,ds
\\
&
\le C_*\delta+C_2\delta^p
\le 2C_*\delta,\qquad 0<t\le 2,
\end{split}
\end{equation}
for any $q\in[r,\infty]$, where $C_1$ and $C_2$ are constants independent of $n_*$ and $\delta$.
Thus we have \eqref{eq:4.20} for $n=n_*+1$, and \eqref{eq:4.20} holds for all $n=1,2,\dots$.

Let $C_{*}^{'}$ be a constant to be chosen later such that $C_{*}^{'}\ge 2C_{*}$. 
By induction we prove that 
\begin{equation}
\label{eq:4.24}
\|u_{n,\lambda}(t)\|_q\le C_{*}^{'}\delta t^{-\frac{N}{2-a}(\frac{1}{r}-\frac{1}{q})},
\quad t>1/2, 
\end{equation}
for any $q\in[r,\infty]$ and $n=1,2,\dots$.
By \eqref{eq:4.19} we have \eqref{eq:4.24} for $n=1$. 
Assume that \eqref{eq:4.22} holds for some $n=n_*$. 
Then, similarly to \eqref{eq:4.23},
since $r_*=\frac{N}{2-a}(p-1)>r$, 
taking a sufficiently small $\delta$ if necessary, 
by (G1), \eqref{eq:4.18} and \eqref{eq:4.20}
we have 
\begin{equation*}
\begin{split}
\|u_{n_*+1,\lambda}(t)\|_q
& \le 
C_3 (t-1/2)^{-\frac{N}{2-a}(\frac{1}{r}-\frac{1}{q})}\|u_{n_*+1,\lambda}(1/2)\|_r
\\
& \qquad\qquad
+C_3\int_{1/2}^{t/2}(t-s)^{-\frac{N}{2-a}(\frac{1}{r}-\frac{1}{q})}\|u_{n_*,\lambda}(s)^p\|_r\,ds
\\
& \qquad\qquad\qquad\qquad
+C_3\int_{t/2}^t\|u_{n_*,\lambda}(s)^p\|_q\,ds
\\
& \le 
C_4C_*\delta t^{-\frac{N}{2-a}(\frac{1}{r}-\frac{1}{q})}
+C_4(C_{*}^{'}\delta)^pt^{-\frac{N}{2-a}(\frac{1}{r}-\frac{1}{q})}
\int_{1/2}^{t/2}s^{-\frac{r_*}{r}}\,ds
\\
& \qquad\qquad\qquad\qquad
+C_4(C_{*}^{'}\delta)^p\int_{t/2}^ts^{-\frac{N}{2-a}(\frac{p}{r}-\frac{1}{q})}\,ds
\\
& \le 
C_5C_*\delta t^{-\frac{N}{2-a}(\frac{1}{r}-\frac{1}{q})}+C_5(C_{*}^{'}\delta)^p 
t^{-\frac{N}{2-a}(\frac{1}{r}-\frac{1}{q})}
\end{split}
\end{equation*}
for all $t>1$, where $C_3$, $C_4$ and $C_5$ are constants independent of $n_*$ and $\delta$.
Let $C_{*}^{'}\ge 2C_5C_{*}$.
Then, taking a sufficiently small $\delta$ if necessary, we have
$$
\|u_{n_*+1,\lambda}(t)\|_q\le C_{*}^{'}\delta t^{-\frac{N}{2-a}(\frac{1}{r}-\frac{1}{q})},\quad t>1. 
$$
This together with \eqref{eq:4.21} implies \eqref{eq:4.24} with $n=n_*+1$.
Thus \eqref{eq:4.24} holds for all $n=1,2,\dots$. 

By \eqref{eq:4.21} and \eqref{eq:4.24} we can find a constant $C$ such that  
$$
\|u_{n,\lambda}(t)\|_q\le C\delta(1+t)^{-\frac{N}{2-a}(\frac{1}{r}-\frac{1}{q})},
\qquad t>0, 
$$
for all $q\in[r,\infty]$ and $n=1,2,\dots$. 
This implies that 
$$
\|u_n(t)\|_q\le C(1+t)^{-\frac{N}{2-a}(\frac{1}{r}-\frac{1}{q})},\qquad t>0, 
$$
for any $q\in[r,\infty]$ and $n=1,2,\dots$.
Then, by the same argument as in the proof of the assertion~(i) of Theorem~\ref{Theorem:1.3}, 
we see that there exists a solution $u$ of \eqref{eq:1.1} satisfying \eqref{eq:1.10}.
Thus the assertion~(ii) of Theorem~\ref{Theorem:1.3} follows,
and the proof of Theorem~\ref{Theorem:1.3} is complete.
$\Box$
\vspace{5pt}

\noindent
{\bf Proof of Corollary~\ref{Corollary:1.1}.}
Since $r_*=N(p-1)/(2-\alpha)$, by \eqref{eq:1.12} 
we can find a constant $C_1$ independent of $\delta$ such that 
$$
\|\varphi\|_{r_*,\infty}\le C_1\delta.
$$
Therefore, by the assertion~(i) of Theorems~\ref{Theorem:1.3} and~\ref{Theorem:1.4} 
we see that,  
if $\delta$ is sufficiently small, 
then a global-in-time solution of \eqref{eq:1.1} exists and it satisfies \eqref{eq:1.8} 
for $\alpha=a$ and \eqref{eq:1.11} for $\alpha=b$.
Thus Corollary~\ref{Corollary:1.1} follows.
$\Box$
\vspace{8pt}

\noindent
{\bf Acknowledgment.}
The first author (Y. F.) was supported in part by the Grant-in-Aid Young Scientists (B) (No. 15K17573)
form Japan Society for the Promotion of Science.
The second author (T. K.) was supported in part by the Grant-in-Aid Young Scientists (B) (No. 16K17629)
form Japan Society for the Promotion of Science.
A part of this work was done while Y. F. and T. K. were visiting the Johns Hopkins University.
They would like to thank the Johns Hopkins University for hospitality.
%


\end{document}